\documentclass[reqno]{amsart}

\usepackage{setspace,enumitem,amsmath,amssymb,amsfonts,color,graphicx,bigints,hyperref}
\usepackage{tikz,array} %This package offers the ability to draw pictures
\usepackage{tikz-cd}

\usepackage[english]{babel}
\usepackage{amsthm}

\usepackage[backend=biber,style=alphabetic, citestyle=alphabetic, sorting = nty]{biblatex}
\theoremstyle{definition}
\newtheorem{theorem}{Theorem}[section]
\newtheorem{lemma}[theorem]{Lemma}
\newtheorem{corollary}[theorem]{Corollary}
\newtheorem{proposition}[theorem]{Proposition}

\theoremstyle{definition}
\newtheorem{definition}[theorem]{Definition}
\newtheorem{remark}[theorem]{Remark}

\newtheorem{conjecture}[theorem]{Conjecture}

\title{Where are the N-Koszul Algebras of Finite Global Dimension?}
\author{Abdourrahmane Kabbaj}
\date{\today}

\begin{document}

\maketitle

\begin{abstract}

The class of $N$-Koszul graded algebras of finite global dimension has gained lots of attention in recent years, especially  in the study of Artin-Schelter regular algebras. While structurally rich and concrete, the only known examples of such algebras are either when $N = 2$, i.e. the algebra is  Koszul, or when $N = 3$. Under a mild Hilbert series assumption, we rule out the existence of $N$-Koszul graded algebras of finite global dimension for $N$ not prime. Furthermore, we establish strong restrictions on the global dimension of such algebras. This suggests that perhaps the existence of 3-Koszul algebras with finite global dimension and  `nice' Hilbert series is an anomaly.

\end{abstract}

\section{Introduction}

Artin-Schelter regular algebras are graded algebras which may be thought of as the noncommutative analogue of polynomial rings. Defined by Artin and Schelter in \cite{Artin1987}, they  are algebras of finite global and GK dimension. Conjectured to have nice  homological and ring theoretic properties  such as being noetherian domains, and having the Hilbert series of a weighted polynomial ring,  their introduction paved the way for the study  of noncommutative projective geometry.  Where noncommutative projective spaces are constructed as noncommutative projective scheme Proj$A$ with $A$ being  an Artin-Schelter regular algebra. Since their introduction, classifying Artin-Schelter regular algebras has been one of most important tasks.  In their original paper,  Artin and Schelter classified  all regular algebras  $A$ of global dimension 3  generated in degree 1 and showed that the minimal projective of the trivial $A$-module $\mathbb{K}$ has one of two shapes:

\[  0 \longrightarrow A( -3) \longrightarrow  A( -2)^3   \longrightarrow  A( -1)^3  \longrightarrow  A  \longrightarrow  \mathbb{K} \longrightarrow 0   \tag{1.0.1}\label{1.0.1} \]

\[  0 \longrightarrow A( -4) \longrightarrow  A( -3)^2   \longrightarrow  A( -1)^2  \longrightarrow  A  \longrightarrow  \mathbb{K} \longrightarrow 0   \tag{1.0.2}\label{1.0.2} \]\\
with those satisfying (\ref{1.0.1}) being Koszul. In an attempt to better understand these algebras, Berger  introduced the following  generalization  of Koszul Algebras that encompasses all the Artin-Schelter regular algebras classified above.\\

\definition\label{def1.1}(\cite{Berger2001}) Let  $A = \bigoplus_{i=1}^\infty A_i$ be an N-graded locally finite graded algebra generated in degree 0 and 1. It is called an \textit{N-Koszul} graded algebra if the trivial $A$-module $\mathbb{K}$ has a minimal projective resolution of the form: 

\[\cdots \longrightarrow P^i \longrightarrow\cdots \longrightarrow  P^1\longrightarrow P^0 \longrightarrow \mathbb{K} \longrightarrow 0\]
for which $P^i$  is generated in degree $\nu(i)$ for all $i\geq 0$ where

 \[
\nu(i) = 
\begin{cases}
    \frac{i}{2}N, & \text{if $i$ is even}  \\
    \frac{i-1}{2}N +1, & \text{if $i$ is odd}
\end{cases}
\]
When $N=2$, this definition coincides with that of $A$ being Koszul.\\

From Berger's definition, we see that all Artin-Schelter regular algebras classified in \cite{Artin1987} are either 2-Koszul, which we refer to as Koszul, or 3-Koszul.  Below is an example of a 3-Koszul Artin-Schelter regular algebra satisfying (\ref{1.0.2}).
 
 \[  A =  \dfrac{\mathbb{K}\langle x,y \rangle }{(xy^2 +y^2x,x^2y + y x^2)}      \]\\

In recent years, the class of $N$-Koszul Artin-Schelter regular algebras has gained significant attention, particularly following Dubois-Violette's demonstration in \cite{DuboisViolette2007} that every $N$-Koszul Artin-Schelter Gorenstein algebra is isomorphic to a derivation quotient algebra $D(\omega, N)$ for some twisted super potential $\omega$, subject to certain pre-regularity conditions. This nice structural  characterization  has lead to an extensive and fruitful study of $N$-Koszul Artin Schelter regular algebras in various different contexts such as \cite{Bocklandt2010}, \cite{Chirvasitu2016}, \cite{Mori2016}, and \cite{Crawford2023}.\\

Although considerable attention has been devoted to the study of these algebras, the only known examples of $N$-Koszul Artin-Schelter regular algebras with $N>2$ remain the 3-Koszul algebras of global dimension 3 first introduced by Artin and Schelter. By contrast, examples of  $N$-Koszul Artin-Schelter Gorenstein algebras  can be found in \cite{DuboisViolette2007}, with the most notable being the infinite family of $N$-Koszul Artin-Schelter Gorenstein algebras  $A(\varepsilon, N)$ for any $N\geq2$ first defined in \cite{Berger2006}. However such examples have been shown to exhibit exponential growth, implying infinite GK dimension for $N>2$ and thus rendering them not Artin-Schelter regular.   Despite their independent significance, it can be deduced from \cite{Stephenson1997} that these algebras lack some of the desirable properties one would like to be true when studying noncommutative projective spaces such as being Noetherian or possessing the Hilbert series of a weighted polynomial ring.\\

The lack of virtually any known examples of $N$-Koszul Artin-Schelter regular thus prompts a question which serves as the motivation for this paper:  Do such algebras really exist  aside from the known 3-Koszul examples?\\

 Let  $R = \mathbb K\langle x_1,x_2,x_3,x_4 \rangle/(f_1,\dots,f_6)$ be the algebra defined in  \cite{Rogalski2012}. While motivated by Artin-Schelter regular algebras, it was shown that $R$ is not  Gorenstein and thus not Artin Schelter regular, yet still posses similar nice properties such as:

\begin{enumerate}

\item[(1)] $R$ is noetherian.

\item[(2)] $R$ is Koszul of finite global dimension $4$

\item[(3)] The Hilbert series of $R$ is $h_R(t) = 1/(1-t)^4$.\\

\end{enumerate}

In an attempt to better understand why examples are so hard to come by, we will drop the Gorenstein condition and more broadly investigate the existence of $N$-Koszul graded algebras of finite global dimension under a mild Hilbert series hypothesis. This paper consists of two mains types of results that we state below. The first establishes restrictions on the possible global dimensions $d$ that can be attained  by $N$-Koszul graded algebra, and the second establishes restrictions on the possible values of $N$ that can appear in examples of $N$-Koszul graded algebra of finite global dimension.\\

 \textbf{For the remainder of the paper, by $\mathbf{N}$-Koszul, we strictly mean $\mathbf{N>2}$.}\\

\theorem\label{theorem1.2} Let $A$ be an $N$-Koszul graded  algebra of finite global dimension $d$ whose Hilbert series is that of a weighted polynomial ring.  If $\mathbb K_A$ has a finite free resolution, then:

\begin{enumerate}

\item[(1)](Proposition \ref{prop3.1} \& Theorem \ref{theorem3.5}) The global dimension is $d =2k+1$ where $k$ is odd.\\

\item[(2)](Corollary \ref{corollary3.9}) The global dimension satisfies  $d \geq  \dfrac{2^N-4}{N} +1$.\\

\item[(3)](Theorem \ref{theorem3.7})  $N$ is prime.\\

\end{enumerate}

Let $A$ be an $N$-Koszul graded algebra of finite global dimension $d$ satisfying the conditions of  the theorems above. The table below illustrates the first handful of possible $N$ values, a lower bound on their respective global dimensions, and whether  an example is known to exist in the literature or not. As seen from the table, the lack of existence of $N$-Koszul algebras of finite global dimension with nice properties might be more than just a mere coincidence.\\

\begin{table}[htbp]
    \centering
\begin{tabular}{|c|c|c|}
    \hline
    
    N & $d \geq$ & Examples Known? \\
    \hline
    2 & 1 & Yes \\
    \hline
    3 & 3 & Yes \\
     \hline
    5 & 7 & No \\
     \hline
    7 & 19 & No \\
     \hline
    11 & 187 & No \\
     \hline
    13& 631 & No \\
     \hline
    17 & 7711 & No \\
     \hline
    19 & 27595 & No \\
     \hline
    \vdots & \vdots & \vdots\\
     \hline
    101 & $2.5\times10^{28}$ & No \\
    \hline
    \end{tabular}
    %\caption{Known examples}
\label{tab:my_table}
\end{table}

These observations lead us to conjecture the following:\\

\conjecture Let $A$ be an $N$-Koszul graded  algebra of finite global dimension $d$ whose Hilbert series is that of a weighted polynomial ring. If  $\mathbb K_A$ has a finite free resolution, then $A$ is a 3-Koszul Artin-Schelter regular algebra of global dimension $3$  classified in \cite{Artin1987}.\\

\section{Preliminaries}

In this section, we review some well known results on graded algebras of finite global dimension, starting with some basic review on graded algebras and modules.\\

Fix a field $\mathbb{K}$.  All vector spaces will be $\mathbb{K}$-vector spaces. An algebra $A$ is called  \textit{$\mathbb N$-graded} if $A = \bigoplus_{n=0}^\infty A_n$ such that $A_iA_j \subseteq A_{i+j}$ for all $i,j$. We say  $A$ is  \textit{locally finite} if $\dim_\mathbb K A_n <\infty$ for all $n$ and \textit{connected} if $A_0 = \mathbb{K}$. Lastly, we say that  $A$ is  \textit{generated in degree one} if all its generators lie in $A_1$. All algebras $A$ in this paper will be assumed to be  $\mathbb N$-graded, locally finite, connected,  and generated  in degree one. \\

A left (respectively right) $A$-module $M$ is \textit{graded} if $M=\bigoplus_{i\in\mathbb{Z}}  M_i$ such that $A_iM_j \subseteq M_{i+j}$ for all $i,j$. For every graded $A$-module $M$ and $n\in \mathbb{N}$, the graded module $M(n) =\bigoplus_{i\in \mathbb{Z}} M(n)_{i}$ where $M(n)_i = M_{n+i}$ for all $i$ denotes the \textit{shift} of $M$ by $n$. A graded $A$-module $M$  is called \textit{left} (respectively \textit{right}) bounded if $M_i= 0$ for $i\ll 0$ (respectively $i \gg0$). We say that $M$ is \textit{bounded} if it is both left and right bounded. For every left bounded graded $A$-module $M$, we define its \textit{Hilbert series} as the formal power series 
\[  h_M(t) = \sum_{i\in \mathbb Z} (\dim_{\mathbb K}M_i)t^i.\]
The Hilbert series of an N-graded algebra $A = \bigoplus_{n=0}^\infty A_n$ is defined similarly.\\

Next  we recall what it means for an  algebra to have finite global dimension. The  \textit{graded projective dimension} of a left $A$-module $M$, denoted by gr.pdim$(M)$  is the minimal (possibly) infinite length of all graded projective resolutions of $M$. The \textit{left global dimension} of $A$, denoted by gr.gldim$_l(A)$, is the supremum of gr.pdim$(M)$ with $M$ ranging over all graded left $A$-modules. The \textit{right global dimension} of $A$, denoted by gr.gldim$_r(A)$ is defined similarly. Since $A$  is connected, it is known that 

\[  \text{gldim}_l(A) =   \text{gr.gldim}_l(A)=\text{pdim}(_A\mathbb{K})  =  \text{pdim}(\mathbb{K}_A)  =  \text{gldim}_r(A) =     \text{gr.gldim}_r(A)   \\     \]

\noindent which shows that one doesn't have to distinguish between the left and right, graded and ungraded global dimensions of $A$ and justifies the use of the terminology ``global dimension of $A$". \\

Suppose that $A$ has finite global dimension $d$. Every graded left bounded projective $A$-module is free, and thus is a sum of shifts of $A$. Thus if $M$ is a left bounded $A$-module, it has a minimal free resolution of the form 

\[  0 \longrightarrow \oplus_{i=1}^{z_d}A(-l_i^d)  \longrightarrow \cdots  \longrightarrow \oplus_{i=1}^{z_1}A(-l_i^1)   \longrightarrow \oplus_{i=1}^{z_0}A(-l_i^0) \longrightarrow M \longrightarrow 0  \tag{2.0.1}\label{2.0.1}\]

\noindent where $z_j$ are possibly infinite. If $z_j$ is finite for all $j$, $M$ is  said to posses a \textit{finite free resolution}. Following the terminology in \cite{Stephenson1997}, if $M$ has a finite free resolution, its \textit{characteristic polynomial} is defined to be 
\[c_M (t) =\sum_{j=0}^d (-1)^j \big( \sum_{i=1}^{z_j} t^{l_i^j} \big).\tag{2.0.2}\label{2.0.2}\]\\

The following lemma connecting the finite  free resolution of an $A$-module $M$ and its Hilbert series will pray a crucial role in our work ahead.\\

\lemma\label{lemma2.1}(\cite{Stephenson1997}, Lemma 2.3) Suppose that $M$ has a finite free resolution of the form (\ref{2.0.1}) and let $c_M(t)$ denote its characteristic polynomial. Then $h_M(t) =c_M(t)h_A(t)$.\\

\remark\label{remark2.2}Let $A$ be an $N$-Koszul graded algebra of finite global dimension $d$ whose trivial $A$-module $\mathbb K_A$ has  a finite free resolution. By Definition \ref{def1.1}, the $\mathbb K_A$ resolution is of the form

\[0 \longrightarrow  \bigoplus_{j=0}^{\beta_d} A(-\nu(d))\longrightarrow \cdots \longrightarrow   \bigoplus_{j=0}^{\beta_2} A(-\nu(2))\longrightarrow   \bigoplus_{j=0}^{\beta_1} A(-\nu(1))\longrightarrow A \longrightarrow \mathbb{K} \longrightarrow 0\]
where $\beta_i$'s are positive integers. Letting $p(t)$ denote the characteristic polynomial of $\mathbb K_A$, it follows from Lemma \ref{lemma2.1}  that $h_{\mathbb K}(t) = p(t)h_A(t)$. Since $h_{\mathbb K} (t) = 1$, it follows that $h_A(t) = \dfrac{1}{p(t)}$ where $p(t) = \sum_{i=0}^d \beta_it^{\nu(i)} $.\\

\definition\label{def2.3} A  graded algebra $A = \bigoplus_{n=0}^\infty A_n$ is said to have the  \textit{Hilbert series of a weighted polynomial ring} if there exists  positive integers $n_i$ such that 
\[  h_A(t) = \frac{1}{ \prod_{i=1}^m(1-t^i)^{n_i}} \]\\

\definition\label{def2.4} A connected graded algebra $A = \bigoplus_{n=0}^\infty A_n$ is called an \textit{ Artin-Schelter regular}  algebra of dimension $d$, AS regular for short, if the following conditions hold:

\begin{itemize}

\item[(1)] $A$ has finite global dimension 

\item[(2)] $A$ is Gorenstein; that is
         \begin{equation*}
    \text{Ext}_A^i(\mathbb{K},A)=
    \begin{cases}
        0 & \text{if } i\neq d\\
        \mathbb K(l) & \text{if } i=d
    \end{cases}
    \end{equation*}
    where $\mathbb K(l)$ is the trivial $A$-module $\mathbb K$ in degree $-l$
            
    \item[(3)] It has finite Gelfand-Kirillov dimension, GK dimension for short; i.e., there are positive numbers $c,d$ such that $\dim A_n \leq cn^d$ for all n.
   
     \end{itemize}       
 An algebra only satisfying conditions (1) and (2) is called an \textit{Artin-Schelter Gorenstein} algebra.\\

Let $A$ be an $N$-Koszul graded algebra of finite global dimension $d$. Throughout our results, we make two assumption  on $A$:\\

\begin{enumerate}

\item[(1)] The Hilbert series of $A$ being that of a weighted polynomial ring.\\

\item[(2)] The trivial $A$-module $\mathbb K_A$ having a finite free resolution.\\

\end{enumerate}

While Remark \ref{remark2.2} and Definition \ref{def1.1} provide a concrete description of the Hilbert series of $A$, the existence of poorly homological and ring theoretic behaved examples such as the family of algebras $A(\varepsilon,N)$ mentioned  in the introduction suggest a need on requiring  $A$ to have  a nice Hilbert series such as that of a weighted polynomial ring. Since this paper is primarily motivated by the class of Artin-Schelter regular algebras defined above. The Hilbert series hypothesis is considered weak as it holds true for all known Artin-Schelter regular algebras and has  been conjectured to be true in general. See \cite{AlexanderPolishchuk2005}, chapter 7.1, remark 3.\\

As to the finiteness of the $\mathbb K_A$ resolution assumption, it is of great importance as most of our results rely on the  ability to express the Hilbert series of $A$ as in Remark \ref{remark2.2}. Notably, this finiteness condition is known to be satisfied for many classes of interest. For instance, 
\begin{enumerate}

\item If $A$ is noetherian, then every finitely generated $A$-module $M$ has a finite  free resolution.

\item If $A$ is Artin-Schelter Gorenstein,   \cite{Stephenson1997} .

\end{enumerate}
and thus both assumptions are considered weak and natural.\\

We end this section by reviewing some results regarding polynomials.\\

Let $f(t)= a_0 + a_1t + \cdots + a_nt^n$ be a polynomial of degree $n$ with integer coefficients. The \textit{conjugate polynomial} of $f(t)$, denoted by $f^*(t)$ is the polynomial $f^*(t) =a_n + a_{n-1}t + \cdots + a_0t^n$. Another equivalent formulation of the conjugate is $f^*(t) = t^nf(t^{-1})$.  A polynomial $f(t)$ is called \textit{self reciprocal} or  \textit{anti-self reciprocal} if $f(t) = f^*(t)$ or $f(t) = -f^*(t)$ respectively. It's easy to see that the  coefficients of a self reciprocal polynomial of degree $n$ satisfy  the condition $a_i =a_{n-i}$ for all $i$. Similarly, the coefficients of an  anti-self reciprocal polynomial of degree $n$ satisfy  the condition $a_i = -a_{n-i}$ for all $i$.\\ 

The lemma below contains some elementary properties of reciprocal polynomials that we state without proof.\\

\lemma\label{lemma2.5} The product of two self reciprocal or anti-self reciprocal polynomials is always self reciprocal. While the product of a self reciprocal polynomial with an anti-self reciprocal polynomial is always anti-self reciprocal.\\

\section{Main Results}

Let $A$ be an $N$-Koszul graded algebra of finite global dimension $d$ whose Hilbert series is that of a weighted polynomial ring. We again emphasize that by $N$-Koszul, we mean $N>2$. In this section, we will establish our main results which are various constraints the pair $(N,d)$ must satisfy.\\

There exist Koszul graded algebras for any global dimension whose Hilbert series is that of a weighted polynomial ring. This next proposition shows this is not the case for  $N$-Koszul graded algebras.\\

\proposition\label{prop3.1} Let $A$ be an $N$-Koszul graded  algebra of finite global dimension $d$ whose Hilbert series is that of a weighted polynomial ring.  If $\mathbb K_A$ has a finite  free resolution, then $A$ is of odd global dimension.

\begin{proof} On the contrary, suppose that  $A$ has even global dimension $d=2k$. It follows from the Hilbert series hypothesis that the characteristic polynomial of $\mathbb K_A$ is $p(t) = \prod_{i=1}^m(1-t^i)^{n_i}$ which by Lemma \ref{lemma2.5} is either self reciprocal or anti-self reciprocal. On the  other hand, 
it follows from  Definition \ref{def1.1} and Remark \ref{remark2.2} that

\[  p(t) =   1 - \beta_1 t + \beta_2 t^N - \beta_3t^{N+1} + \cdots - \beta_{d-1} t^{(k-1)N+1} + \beta_dt^{kN}   \tag{3.1.1}\label{3.1.1}     \]\\
where  $\beta_i$'s are positive integers. Since $p(t)$ being anti-self reciprocal would imply that $\beta_d = -1$ which isn't possible, we conclude that it's self reciprocal and thus $p(t) = t^{kN}p(t^{-1})$. Expanding both sides and comparing the exponents of the first non constant term  yields that $N=2$ which is a contradiction.

\end{proof}

\remark\label{remark3.2} Let $A$ be as in Theorem \ref{theorem1.2}. By the previous proposition, $A$ is of odd global dimension $d=2k+1$ for some positive integer $k$.  We will frequently express the characteristic polynomial  $p(t)$ of the trivial $A$-module $\mathbb K_A$ in one of the following three equivalent forms:

\begin{itemize}

\item[(1)] By the Hilbert series hypothesis,  
 \[p(t) = \prod_{i=1}^m(1-t^i)^{n_i} \tag{3.2.1}\label{3.2.1}\]

\item[(2)] By Remark \ref{remark2.2} and the proposition above,  
\[ p(t) = 1 - \beta_1 t + \beta_2 t^N - \beta_3 t^{N+1} +\;\; \cdots  \;\;- \beta_2t^{(k-1)N+1} + \beta_1 t^{kN} - t^{kN+1}  \tag{3.2.2}\label{3.2.2}\]

\item[(3)] Letting  $q(t) = \sum_{i=0}^k \beta_{2i}t^{iN}$, the formulation in (2) can we rewritten as  

\[  p(t) = q(t) - tq^*(t)    \tag{3.2.3}\label{3.2.3}\]
where $q^*(t)$ denotes the conjugate of $q(t)$.\\

\end{itemize}

We observe from the remark above that the  Hilbert series hypothesis necessitates that all the roots of $p(t)$ must be roots of unity. Combining this with the fact that those same roots must satisfy the equivalent formulations in (\ref{3.2.2}) and (\ref{3.2.3}), one can establish some restrictions on the roots of $p(t)$ which will come in handy when proving some of our main theorems.\\

\proposition\label{prop3.3} Let $A$ be an $N$-Koszul graded  algebra of finite global dimension $d$ whose Hilbert series is that of a weighted polynomial ring.  If $\mathbb K_A$ has a finite  free resolution and  $\zeta_j$ is  a primitive $j^{th}$ root of unity which is also a root of $p(t)$, then $\gcd(j,N) = 1$.\\

\begin{proof} On the contrary, suppose that $\zeta_j$ is a primitive  root of unity which is also a root of $p(t)$ yet $\gcd(j, N) = s>1$. By  the Hilbert series assumption, this implies that $1- t^j$ divides $p(t)$. Since $s$ divides $j$, $ 1 - t^s$ also divides $p(t)$, implying that $\zeta_s$ is also a root of $p(t)$. Since $s$  divides $N$, it follows that $\zeta_s^N = 1$ and $q(\zeta_s) = q^*(\zeta_s) = \sum_{i=0}^k \beta_{2i}>0$ since  all the $\beta_{2i}$ are positive integers. By (\ref{3.2.3}),
 \[ p(\zeta_s) =0 \implies q(\zeta_s) = \zeta_s q^*(\zeta_s)\]
 
 \[ \implies \zeta_s = \frac{q(\zeta_s)}{q^*(\zeta_s)}   =   1 \]
 which is a contradiction.  
 \end{proof}

\lemma\label{lemma3.4} Let $A$ be an $N$-Koszul graded  algebra of finite global dimension $d$ whose Hilbert series is that of a weighted polynomial ring.  If $\mathbb K_A$ has a finite free resolution and  $\zeta_j$ is  a primitive $j^{th}$ root of unity which is also a root of $p(t)$ written as in (\ref{3.2.3}), then $q(\zeta_j)\neq 0$.
 
\begin{proof}  On the contrary, suppose that $\zeta_j$ is a primitive  root of unity which is also a root of $p(t)$ expressed as in (\ref{3.2.3}) yet $q(\zeta_j)=0$. Then $q^*(\zeta_j)=0$ as well. Define $\alpha = \zeta_j\cdot \zeta_N$. Since \[ \alpha^{Ni} = (\zeta_j^{Ni})\cdot(\zeta_N)^{Ni} =\zeta_j^{Ni}\;\;\;\;\;\; \text{for all $i$ }    \]
 its easy to see that $q(\alpha) = q(\zeta_j) = q^*(\alpha) = q^*(\zeta_j) =0$, and so $p(\alpha)=0$. Since $\gcd(j,N) =1$ by Proposition \ref{prop3.3},  $\alpha$ is a primitive $jN^{th}$ root of unity and is also a root of $p(t)$ but this is a contradiction by Proposition \ref{prop3.3}  since $\gcd(jN,N) \neq 1$. 
 
 \end{proof}

In \cite{DuboisViolette2007}, it was shown  that if an algebra $A$ is $N$-Koszul Artin-Schelter Gorenstein  of finite  global dimension, then $d=2k+1$. Therefore, Proposition \ref{prop3.1} can be seen as an extension of that result to a broader context. This next theorem shows that one can actually establish more restrictions on $d$.\\

\theorem\label{theorem3.5} Let $A$ be an $N$-Koszul graded  algebra of finite global dimension $d$ whose Hilbert series is that of a weighted polynomial ring.  If $\mathbb K_A$ has a finite free resolution, then  $d =4r+3$ for some positive integer $r$.\\

\begin{proof}  By our Hilbert series assumption, 
\[    p(t) = (1-t)^{n_1}(1-t^2)^{n_2}(1-t^3)^{n_3}\prod_{i=4}^m (1-t^{i})^{n_i}   \]

\[   = 1 - n_1t + \Big[{n_1\choose 2} - n_2 \Big]t^2 + \Big[ -{n_1\choose 3} + n_1n_2 -n_3\Big]t^3 + \cdots  \tag{3.5.1}\label{3.5.1}\]\\

From (\ref{3.5.1}) and the relationship between $p(t)$ and the $\mathbb{K}_A$ projective resolution of $A$, we see  that $n_1$ corresponds to the number of generators of our algebra $A$. Since $A$ is clearly not finite dimensional as a $\mathbb{K}$-vector space, we can  assume that $n_1\ge 2$ as otherwise $A\cong \mathbb K[x]$ which is known to be Koszul. By comparing (\ref{3.5.1})   with (\ref{3.2.2}) and using the fact that $N>2$, we conclude that ${n_1\choose 2} - n_2 =0$. Thus $n_2 = {n_1\choose 2} \geq 1$ which then implies that $-1$ is a root of $p(t)$.  It then follows by (\ref{3.2.3}) that 

\[p(-1) = 0 \implies  q(-1) + (- 1)^{kN}q(-1) = q(-1) + q^*(-1)  = 0\]
\[\implies q(-1)[1 + (-1)^{kN}] =0\]
Since  $q(-1) \neq 0$  by Lemma \ref{lemma3.4}, it follows that $kN$ is odd and thus both $k$ and $N$  are odd as well. It then follows that $d=2k+1 = 2(2r+1 ) = 4r+3$ for some positive integer $r$.

\end{proof}

One can actually say a bit more from the theorem above. If we recall the family $A(\varepsilon, N)$ of $N$-Koszul Artin-Schelter Gorenstein algebras mentioned in the introduction, it was shown that examples exist for any $N\geq2$. In the setting of our theorem above, restrictions are now being imposed on the possible $N$ values. For instance, one shouldn't expect to find an example of a 4-Koszul Artin-Schelter regular algebra. We will see from the next two results that we can further establish much stronger restrictions on the possible values on $N$.\\

For this next lemma, we will use the following Taylor series expansion: 

$$    -\log(1-t) = \sum_{j=1}^\infty \dfrac{t^j}{j}   $$

\lemma\label{lemma3.6} Suppose that \[ \prod_{i=1}^m(1-t^i)^{n_i} =1 - n_1 t + \beta_2 t^N - \beta_3 t^{N+1} +\;\; \cdots  \;\;- t^{kN+1} \tag{3.6.1}\label{3.6.1}\]  
If $n_1\geq2$,  then $n_r \geq \dfrac{n_1^{r-1}}{r}$ for $1\leq r\leq N-1$.\\

\begin{proof}

 Apply $-\log$ to both sides of $(\ref{3.6.1})$ to obtain the following: 

\[  -\log\Big(\prod_{i=1}^m (1-t^i)^{n_i}  \Big)  =     -\log\Big( 1 - n_1t + \beta_2t^N - \beta_3t^{N+1} +\cdots -t^{kN +1} \Big) \tag{3.6.2}\label{3.6.2}     \]
Expanding the left hand side of (\ref{3.6.2}), we get 

\[ -\log\Big(\prod_{i=1}^m (1-t^i)^{n_i}  \Big)   =\sum_{i=1}^m - n_i\log(1-t^i)             \]

\[ = \sum_{i=1}^m n_i\sum_{j=1}^\infty  \dfrac{t^{ij}}{j}  \tag{3.6.3}\label{3.6.3}   \]
For each $1 \leq r \leq N-1$, it follows from (\ref{3.6.3}) that the coefficient of $t^r$  on the left hand side of $(\ref{3.6.2})$ is  \[ \sum_{i|r}   n_i \dfrac{1}{\Big(\dfrac{r}{i}\Big)} = \dfrac{1}{r} \sum_{i|r} in_i               \]\\
Now to expand the  right hand side of $(\ref{3.6.2})$, let $u = n_1t - \beta_2t^N +\beta_3t^{N+1} - \cdots + t^{kN+1}$. Then 

\[   -\log\Big( 1 - n_1t + \beta_2t^N - \beta_3t^{N+1} +\cdots -t^{kN +1} \Big)   = -\log(1-u)    \]

\[ = \sum_{j=1}^\infty \dfrac{u^j}{j}  =  \sum_{j=1}^\infty \dfrac{( n_1t - \beta_2t^N +\beta_3t^{N+1} - \cdots + t^{kN+1})^j}{j}     \]

\[   = \sum_{j=1}^\infty \dfrac{(n_1t)^j}{j} + t^N\sum_{l=0}^\infty \alpha_lt^l     \tag{3.6.4}\label{3.6.4} \]
for some $\alpha_l$. For each $1 \leq r \leq N-1$, it follows from (\ref{3.6.4}) that the coefficient of $t^r$  on the right hand side of (\ref{3.6.2}) is $\dfrac{n_1^r}{r}$. Since the coefficients  have to match on  both sides of (\ref{3.6.2}), we obtain the following recursive formula:

\[   n_1^r=  \sum_{i|r}in_i  \tag{3.6.5}\label{3.6.5} \]
for all $1 \leq r\leq N-1$. Applying the Möbius inversion formula to (\ref{3.6.5}), we get 

\[  rn_r = \sum_{i|r} \mu(i)n_1^{\frac{r}{i}}    \]
By  hypothesis and the discussion in Theorem \ref{theorem3.5}, the inequality is satisfied for $r=1$ and $r=2$. For $3 \leq r \leq N-1$, let $p$ be a minimal prime such that $p |r$. Then 

\[  rn_r = \sum_{i|r} \mu(i)n_1^{\frac{r}{i}}  = n_1^r +   \sum_{1< i|r} \mu(i)n_1^{\frac{r}{i}} \]

\[ \geq   n_1^r -   \sum_{1< i|r} n_1^{\frac{r}{i}}  = n_1^r -    \sum_{\substack{ 1\leq j\leq \frac{r}{p},\\j |r}} n_1^j \]

\[  \geq    n_1^r -    \sum_{j=0}^ \frac{r}{p}n_1^j  \; = \;  n_1^r  - \dfrac{n_1^ {\frac{r}{p}+1} -1}{ n_1 -1}\]

\[  \geq    n_1^r - n_1^ {\frac{r}{p}+1} +1 \geq n_1^r - n_1^{r-1}  \geq n_1^{r-1}  \]

\[ \implies n_r \geq \frac{n_1^{r-1}}{r}\qedhere \] 
\end{proof}

 We are now finally ready  to prove one of our main results establishing strong restrictions on $N$.
 
\theorem\label{theorem3.7} Let $A$ be an $N$-Koszul graded  algebra of finite global dimension $d$ whose Hilbert series is that of a weighted polynomial ring.  If $\mathbb K_A$ has a finite free resolution, then $N$ is prime.\\
 
\begin{proof} On the contrary, suppose that $N =rs$ is a composite number. By the discussion in Remark \ref{remark3.2}, the characteristic polynomial of $\mathbb{K_A}$ is
 \[p(t) =  \prod_{i=1}^m(1-t^i)^{n_i} =1 - \beta_1 t + \beta_2 t^N - \beta_3 t^{N+1} +\;\; \cdots  \;\;- t^{kN+1} \]  
 Comparing the first nonconstant  term in both sides yields that $\beta_1 = n_1$ which is greater than or equal to two by similar reasoning as Theorem \ref{theorem3.5}. By Lemma \ref{lemma3.6}, it follows that   $n_1,\cdots, n_{N-1} >0$ from which we conclude that the primitive  $j^{th}$ root of unity $\zeta_j$ is a  root of $p(t)$ for all $1\leq j\leq N-1$. In particular, $\zeta_r$ is a primitive root of $p(t)$.  However, this is a  contradiction since $\gcd(r,N) \neq1$ which isn't possible by Proposition \ref{prop3.3}. So we conclude that N is prime.\\

\end{proof}

Going back to our study of the  global dimensions  of our $N$-Koszul algebras,  the next theorem establishes a lower bound for the global dimension which is exponential in $N$.\\

\theorem\label{theorem3.8} Let $A$ be an $N$-Koszul graded  algebra of finite global dimension $d$ whose Hilbert series is that of a weighted polynomial ring.  If $\mathbb K_A$ has a finite free resolution, then  $d  \geq  \dfrac{2 (s^{N-2}-1)}{N} +1 $ where $s$ denotes the number of generators of $A$.
 
\begin{proof} By the discussion in Remark \ref{remark3.2}, the characteristic polynomial of $\mathbb{K_A}$ is  \[p(t) =  \prod_{i=1}^m(1-t^i)^{n_i} =1 - \beta_1 t + \beta_2 t^N - \beta_3 t^{N+1} +\;\; \cdots  \;\;- t^{kN+1} \]  
 Since $A$ is generated by $s\geq 2$ elements, Lemma \ref{lemma3.6} applied here yields that   $n_i \geq \dfrac{s^{i-1}}{i}$ for $1\leq i \leq N-1$ and thus $f(t) = \prod_{i=1}^{N-1}(1-t^i)^{n_i}$  must divide  $p(t)$.  This implies that $\deg f(t) \leq \deg p(t)$ and so

\[\dfrac{s^{N-1}  -1}{s-1} =  \sum_{i=1}^{N-1} i\Big(\dfrac{s^{i-1}}{i}\Big)  \leq \sum_{i=1}^{N-1}in_i   = \deg f(t) \leq \deg p(t)  =  kN +1\]

\[\implies \dfrac{d-1}{2} = k \geq  \dfrac{s\big(s^{N-2}-1\big)}{(s-1)N} \geq   \dfrac{s^{N-2}-1}{N} \]

\[\implies d  \geq  \dfrac{2 (s^{N-2}-1)}{N} +1\\ \qedhere \]

\end{proof}

\corollary\label{corollary3.9} Let $A$ be an $N$-Koszul graded  algebra of finite global dimension $d$ whose Hilbert series is that of a weighted polynomial ring.  If $\mathbb K_A$ has a finite free resolution, then  $d  \geq  \dfrac{2^N-4}{N} +1 $.

\begin{proof} Since every $N$-Koszul graded  algebra is generated by at least two elements. It follows from the proof of Theorem \ref{theorem3.8} that 
\[   2^{N-1} -1 \leq \sum_{i=1}^{N-1} in_i \leq kN+1        \]
\[ \implies d  \geq  \dfrac{2^N-4}{N} +1 \qedhere \]

\end{proof}

With the only known example of $N$-Koszul Artin-Schelter regular algebras being 3-Koszul of global dimension 3. One might expect to find other $d$-Koszul graded algebras of global dimension $d$. This next corollary establishes otherwise.\\

\corollary Let $A$ be an $d$-Koszul graded  algebra of finite global dimension $d$ whose Hilbert series is that of a weighted polynomial ring.  If $\mathbb  K_A$ has a finite free resolution, then $d=3$.
 
\begin{proof} This follows immediately from the corollary above  since $d \geq   \dfrac{2^d-4}{d} +1 > d$ for $d$ greater than three.

\end{proof}

\remark Suppose one is now interested in constructing new examples of $N$-Koszul graded algebras of finite global dimension, perhaps in the study of Artin-Schelter regular algebras. An idea could be to construct such examples as some twisted tensor product of two known $N$-Koszul graded algebras $A$ and $B$. This is easily seen to not work as any form of a twisted tensor product of two algebras generated in degree one will introduce quadratic relations which cannot occur in an $N$-Koszul algebra.\\

It seems likely that an attempt for constructing new examples, if they exist,  would require explicit projective resolution constructions which is known to a very challenging task for algebras of large global dimension which happens to be the case.   \\

%This subsequent result illustrates that such constructions wont work and thus one might need to directly compute an explicit representation which is a very challenging task especially in large global dimensions. For more on twisted tensor products, we refer the reader to \cite{Shepler2016}.

%\proposition Let $A$  and $B$ be  N$_1$-Koszul graded algebras of global dimension $d_1$ and  N$_2$-Koszul graded algebras of global dimension $d_2$ respectively. For all strictly graded twisting maps $\tau$, the algebra  $A \otimes_\tau B$ is not $N$-Koszul.\\

%\begin{proof} By Remark \ref{remark3.2}, there exists positive integers $k_1$ and $k_2$ such that $d_1=2k_1 +1$ and $d_2 =2k_2 +1$. Furthermore,  there exists positive integers $\alpha_1,\cdots,\alpha_{d_1}$ and $\beta_1,\cdots,\beta_{d_2}$ such that $h_A(t) = \frac{1}{p_A(t)}$ and $h_B(t) = \frac{1}{p_B(t)}$ where 

%\[  p_A(t) =   1 - \alpha_1 t + \alpha_2 t^N - \cdots - \alpha_{d_1-1} %t^{(k_1-1)N+1} + \alpha_{d_1}t^{k_1N}       \]

%\[  p_B(t) =   1 - \beta_1 t + \beta_2 t^N - \cdots - \beta_{d_2-1} t^{(k_2-1)N+1} + \beta_{d_2}t^{k_2N}       \]

%Since $\tau$ is strongly graded, it follows from \cite{Shepler2016},  and the \text{Künneth}  formula that

%\[   h_{A\otimes_\tau B}(t) = h_{A\otimes_\mathbb K B}(t) =  h_A(t)h_B(t)   \]

%and thus  
%\[p_{A\otimes_{\tau} B}(t) = p_A(t)\cdot p_B(t)\]
%\[   =1 - (\alpha_1 + \beta_1)t + \alpha_1\beta_1 t^2 + \cdots +      %\alpha_{d_1}\beta_{d_2}t^{(k_1+k_2)N}    \]

%which cannot correspond to an $N$-Koszul graded algebra as this %characteristic polynomial has a second degree term. 

% \end{proof}

 \section*{Acknowledgments}

The author thanks Vladimir Baranovsky and So Nakamura for helpful conversations that lead to a proof of Lemma 3.6. In addition, the author would like to thank his advisor Manny Reyes for his invaluable mentorship, feedback, and proofreading of this paper. This work constitutes a portion of the authors Ph.D. thesis.

\printbibliography

@Article{Artin1987,
  author  = {Michael Artin and William F Schelter},
  journal = {Advances in Mathematics},
  title   = {Graded algebras of global dimension 3},
  year    = {1987},
  issn    = {0001-8708},
  number  = {2},
  pages   = {171-216},
  volume  = {66},
  doi     = {https://doi.org/10.1016/0001-8708(87)90034-X},
  url     = {https://www.sciencedirect.com/science/article/pii/000187088790034X},
}

@Article{Chirvasitu2016,
  author  = {Alexandru Chirvasitu and Chelsea M. Walton and Xingting Wang},
  journal = {Journal of Noncommutative Geometry},
  title   = {On quantum groups associated to a pair of preregular forms},
  year    = {2016},
  url     = {https://api.semanticscholar.org/CorpusID:119704612},
}

@Article{Bocklandt2010,
  author   = {Raf Bocklandt and Travis Schedler and Michael Wemyss},
  journal  = {Journal of Pure and Applied Algebra},
  title    = {Superpotentials and higher order derivations},
  year     = {2010},
  issn     = {0022-4049},
  number   = {9},
  pages    = {1501-1522},
  volume   = {214},
  doi      = {https://doi.org/10.1016/j.jpaa.2009.07.013},
  url      = {https://www.sciencedirect.com/science/article/pii/S002240490900173X},
}

@Article{DuboisViolette2007,
  author   = {Michel Dubois-Violette},
  journal  = {Journal of Algebra},
  title    = {Multilinear forms and graded algebras},
  year     = {2007},
  issn     = {0021-8693},
  number   = {1},
  pages    = {198-225},
  volume   = {317},
  doi      = {https://doi.org/10.1016/j.jalgebra.2007.02.007},
  url      = {https://www.sciencedirect.com/science/article/pii/S0021869307000968},
}

@Article{Mori2016,
  author   = {Izuru Mori and S. Paul Smith},
  journal  = {Journal of Algebra},
  title    = {m-Koszul Artin–Schelter regular algebras},
  year     = {2016},
  issn     = {0021-8693},
  pages    = {373-399},
  volume   = {446},
  doi      = {https://doi.org/10.1016/j.jalgebra.2015.09.016},
  url      = {https://www.sciencedirect.com/science/article/pii/S0021869315004615},
}

@Book{AlexanderPolishchuk2005,
  author    = {Alexander Polishchuk, Leonid Positselski},
  publisher = {American Mathematical Soc.},
  title     = {Quadratic Algebras},
  year      = {2005},
  isbn      = {0821838342},
}

@Article{Berger2006,
  author   = {Berger, Roland and Marconnet, Nicolas},
  journal  = {Algebras and Representation Theory},
  title    = {Koszul and Gorenstein Properties for Homogeneous Algebras},
  year     = {2006},
  issn     = {1572-9079},
  number   = {1},
  pages    = {67--97},
  volume   = {9},
  doi      = {10.1007/s10468-005-9002-1},
  refid    = {Berger2006},
  url      = {https://doi.org/10.1007/s10468-005-9002-1},
}

@Article{Rogalski2012,
  author  = {Rogalski, D. and Sierra, Susan J.},
  journal = {Compositio Mathematica},
  title   = {Some projective surfaces of GK-dimension 4},
  year    = {2012},
  number  = {4},
  pages   = {1195–1237},
  volume  = {148},
  doi     = {10.1112/S0010437X12000188},
}

@Article{Berger2001,
  author  = {Roland Berger},
  journal = {Journal of Algebra},
  title   = {Koszulity for Nonquadratic Algebras},
  year    = {2001},
  pages   = {705-734},
  volume  = {239},
  url     = {https://api.semanticscholar.org/CorpusID:120370197},
}

@InProceedings{Stephenson1997,
  author = {Darin R. Stephenson and James J. Zhang},
  title  = {Growth of graded noetherian rings},
  year   = {1997},
  url    = {https://api.semanticscholar.org/CorpusID:15994131},
}

@Article{Crawford2023,
  author   = {Crawford, Simon},
  journal  = {Algebras and Representation Theory},
  title    = {Superpotentials and Quiver Algebras for Semisimple Hopf Actions},
  year     = {2023},
  issn     = {1572-9079},
  number   = {6},
  pages    = {2709--2752},
  volume   = {26},
  doi      = {10.1007/s10468-022-10165-y},
  refid    = {Crawford2023},
  url      = {https://doi.org/10.1007/s10468-022-10165-y},
}

\end{document}